# A NEW CLASS OF GENERALIZED BAYES MINIMAX RIDGE REGRESSION ESTIMATORS

BY YUZO MARUYAMA AND WILLIAM E. STRAWDERMAN[1]

*University of Tokyo and Rutgers University*

Let $y = A\beta + \varepsilon$, where $y$ is an $N \times 1$ vector of observations, $\beta$ is a $p \times 1$ vector of unknown regression coefficients, $A$ is an $N \times p$ design matrix and $\varepsilon$ is a spherically symmetric error term with unknown scale parameter $\sigma$. We consider estimation of $\beta$ under general quadratic loss functions, and, in particular, extend the work of Strawderman [*J. Amer. Statist. Assoc.* **73** (1978) 623–627] and Casella [*Ann. Statist.* **8** (1980) 1036–1056, *J. Amer. Statist. Assoc.* **80** (1985) 753–758] by finding adaptive minimax estimators (which are, under the normality assumption, also generalized Bayes) of $\beta$, which have greater numerical stability (i.e., smaller condition number) than the usual least squares estimator. In particular, we give a subclass of such estimators which, surprisingly, has a very simple form. We also show that under certain conditions the generalized Bayes minimax estimators in the normal case are also generalized Bayes and minimax in the general case of spherically symmetric errors.

**1. Introduction.** In this paper we consider adaptive ridge regression estimators in the general linear model with homogeneous spherically symmetric errors. There are three main contributions: (a) we propose sufficient conditions on estimators for simultaneously reducing risk and increasing numerical stability relative to the least squares estimator for all full rank design matrices, (b) under normality, we obtain a broad class of generalized Bayes estimators satisfying the above sufficient conditions, and (c) this class contain a subclass of particularly simple form, which, we hope, adds to the practical utility of our results.

Hoerl and Kennard [11] introduced the ridge regression technique as a way to simultaneously reduce the risk and increase the numerical stability of the

Received December 2003; revised September 2004.
[1]Supported in part by NSA Grant H98230-04-0114.

least squares estimator in ill-conditional problems. The risk reduction aspect of Hoerl and Kennard's method was often observed in simulations but was not theoretically justified. Strawderman [19] looked at the problem in the context of minimaxity and produced minimax adaptive ridge-type estimators, but ignored the condition number aspect of the problem. Casella [7, 8] considered both the minimaxity and condition number aspects and gave estimators which were minimax and condition number decreasing for some, but not all, design matrices. Neither Strawderman nor Casella gave generalized Bayes minimax estimators. Moreover, to the best of our knowledge, almost all theoretical results on ridge regression in the literature depend on normality.

In the present paper we propose a broad class of minimax estimators which increases the numerical stability of the least squares estimator for all full rank design matrices, under the assumption of a spherically symmetric error distribution. Furthermore, under normality a broad class of generalized Bayes estimators included in the above class is found. What is particularly noteworthy about our class of estimators is that it contains a subclass with a form (adapted to the case of unknown $\sigma^2$) which is remarkably similar to that of the estimators originally suggested in [17] for the case $\mathrm{Cov}(X) = I$. In particular, our simple generalized Bayes estimators of the mean vector are of the form

$$\hat{\theta}_{\mathrm{SB}} = (I - \alpha/\{\gamma(\alpha+1) + W\}C^{-1})X,$$

where $W = X'C^{-1}D^{-1}X/S$ for some positive-definite matrices $C$ and $D$.

To be more precise, we start the familiar linear regression model $Y = A\beta + \varepsilon$, where $Y$ is an $N \times 1$ vector of observations, $A$ is the known $N \times p$ design matrix of rank $p$, $\beta$ is the $p \times 1$ vector of unknown regression coefficients, and $\varepsilon$ is an $N \times 1$ vector of experimental errors. We assume $\varepsilon$ has a spherically symmetric distribution with a density $\sigma^{-N}f(\varepsilon'\varepsilon/\sigma^2)$, where $\sigma$ is an unknown scale parameter and $f(\cdot)$ is a nonnegative function on the nonnegative real line.

The least squares estimator of $\beta$ is $\hat{\beta} = (A'A)^{-1}A'y$. Since the covariance matrix of $\hat{\beta}$ is proportional to $\sigma^2(A'A)^{-1}$, the least squares estimator may not be a suitable estimator when some components of $\hat{\beta}$ or some linear combinations of $\hat{\beta}$ have a very large variance and when $A'A$ is nearly singular. Additionally, $(A'A)^{-1}$ may have inflated diagonal values so that small changes in the observations produce large changes in $\hat{\beta}$. Hoerl and Kennard [11] proposed the ridge estimator

$$(1.1) \qquad \hat{\beta}_{\mathrm{R}}(k) = (A'A + kI_p)^{-1}A'y,$$

where $k$ is a positive constant, to ameliorate these problems. Adding the number $k$ before inverting amounts to increasing each eigenvalue of $A'A$ by $k$.

In particular, if $P$ is the $p \times p$ orthogonal matrix of eigenvectors of $(A'A)^{-1}$, with $d_1 \geq d_2 \geq \cdots \geq d_p$ as eigenvalues, it follows that

$$P'(A'A)^{-1}P = D, \qquad P'P = I_p,$$

where $D = \mathrm{diag}(d_1, \ldots, d_p)$. Then (1.1) can be written as

(1.2) $$\hat{\beta}_{\mathrm{R}}(k) = P(D^{-1} + kI_p)^{-1}P'Ay.$$

The ridge estimator is more stable than $\hat{\beta}$ in the sense that the condition number of the estimator is reduced.

However, we are interested in proposing better estimators than $\hat{\beta}$ from the decision-theoretic point of view. We measure the loss in estimating $\beta$ by $b$ with loss functions

(1.3) $$L_j(b, \beta, \sigma^2) = \sigma^{-2}(b - \beta)'(A'A)^j(b - \beta),$$

where $(A'A)^j = P\,\mathrm{diag}(d_1^{-j}, \ldots, d_p^{-j})P'$. In particular, $L_j$ for $j = 0, 1, 2$ are known as squared error loss, predictive (or scale invariant) loss and Strawderman's [19] loss, respectively. Then the risk function of an estimator $b$ is given by $R_j(b, \beta, \sigma^2) = E[L_j(b, \beta, \sigma^2)]$. The least squares estimator $\hat{\beta}$ is minimax with constant risk. Therefore, $b$ is a minimax estimator of $\beta$ if and only if $R_j(b, \beta, \sigma^2) \leq R_j(\hat{\beta}, \beta, \sigma^2)$ for all $\beta$ and $\sigma^2$. Hence, the search for estimators better than $\hat{\beta}$ is a search for minimax estimators.

To simplify expressions and to make matters a bit clearer, it is helpful to rotate the problem via the following transformation, so that the covariance matrix of $\hat{\beta}$ becomes diagonal. Let $Q$ be an $N \times N$ orthogonal matrix such that

$$QA = \begin{pmatrix} D^{-1/2}P' \\ 0 \end{pmatrix}$$

and let $D_*$ be the $N \times N$ diagonal matrix $\mathrm{diag}(d_1, \ldots, d_p, 1, \ldots, 1)$. Next define two random vectors $X = (X_1, \ldots, X_p)'$ and $Z = (Z_1, \ldots, Z_n)'$, where $n = N - p$, by

$$\begin{pmatrix} X \\ Z \end{pmatrix} = D_*^{1/2}QY.$$

Then $(X', Z')'$ has the joint density given by

(1.4) $$\prod d_i^{-1/2} \sigma^{-p-n} f(\{(x - \theta)'D^{-1}(x - \theta) + z'z\}/\sigma^2),$$

where $\theta = P'\beta$. Notice also that $X$ and $Z'Z$ can be expressed as $P'\hat{\beta}$ and $(y - A\hat{\beta})'(y - A\hat{\beta})$, respectively. Denote $Z'Z$ by $S$, as is customary. The original problem is thus equivalent to estimation of $\theta$ under the loss function $L_j(\delta, \theta, \sigma^2) = (\delta - \theta)'D^{-j}(\delta - \theta)/\sigma^2$, where $j$ corresponds to $j$ in (1.3). We will consider the problem in this equivalent canonical form. Note that $L_0$ is, in a sense, the least favorable among $L_j$ for $j \geq 0$ under multicollinearity because $L_j$ for $j > 0$ relatively reduces the contribution of components with large variance.

Strawderman [19] and Casella [7] essentially considered the class of estimators of the form

$$\hat{\theta}_{\mathrm{R}}(K) = (I - \{I + D^{-1}K^{-1}\}^{-1})X,$$

which originally came from straight generalization of (1.2), that is, the generalized ridge estimator

(1.5) $\qquad \hat{\beta}_{\mathrm{R}}(K) = (A'A + PKP')^{-1}A'y = P(D^{-1} + K)^{-1}P'A'y,$

where $K = \mathrm{diag}(k_1, \ldots, k_p)$. Under general quadratic loss they proposed a sufficient condition for minimaxity under normality for adaptive estimators $\hat{\theta}_{\mathrm{R}}(\hat{K})$, where $\hat{K} = \psi(X'D^{-1}X/S)\mathrm{diag}(a_1, \ldots, a_p)$, $\psi$ is a suitable positive function and $a_i$ is positive for all $i$. Casella [7] discussed the relationship between minimaxity and stability (in terms of lowered condition number) and pointed out that forcing ridge regression estimators to be minimax makes it difficult for them to provide the numerical stability for which they were originally intended. Casella [8] found that, under certain conditions on the structure of the eigenvalues of the design matrix, both minimaxity under $L_0$ and stability can be simultaneously achieved for a special case $\psi(w) \equiv w^{-1}$.

In Section 2, for the general spherically symmetric case, we give a class of minimax estimators of $\theta$ (and hence, by transformation, $\beta$) under $L_j$, somewhat broader than those of Strawderman [19] and Casella [7, 8]. We then give a class of generalized hierarchical prior distributions on $\theta$ and $\sigma^2$ which, in the normal case, give generalized Bayes estimators satisfying the minimaxity condition. This class generalizes (also to the class of unknown $\sigma^2$) the class of priors in [3, 4, 10, 14, 18]. We further show that, for certain choices of parameters in the hierarchy, the resulting estimators have the simple form indicated above. We also show that in certain cases a version of our minimax estimator is generalized Bayes for the entire class of spherically symmetric error distributions. Section 3 is devoted to the study of general conditions under which the generalized ridge regression estimator $\hat{\beta}_{\mathrm{R}}(K)$ competitive with $\hat{\beta}$ has increased numerical stability (i.e., decreased condition number). Section 4 is devoted to showing that we may always choose a minimax estimator (which is also generalized Bayes under normality) in our class which has greater numerical stability than the least squares estimator. In particular, our simple generalized Bayes minimax stable estimators under normality are quite practical for the general spherically symmetric case. In Section 5 we give some numerical results.

**2. A class of minimax generalized Bayes estimators.** In this section we first give a sufficient condition for minimaxity under the loss $L_j$ and the spherically symmetric case, and then use it to obtain a class of generalized Bayes minimax estimators under the normal case. This class contains a subclass of a particularly simple form, which we hope adds to the practical utility of our results. We also show that in certain cases a version of our minimax estimator is generalized Bayes for the entire class of spherically symmetric error distributions.

Our estimators are of the form

$$\hat{\theta}_\phi = \left(I - \frac{S}{X'C^{-1}D^{-1}X}\phi\left(\frac{X'C^{-1}D^{-1}X}{S}\right)C^{-1}\right)X, \tag{2.1}$$

where $C = \text{diag}(c_1, \ldots, c_p)$, where $c_i \geq 1$ for any $i$. We note that estimators of the form (2.1) satisfy "directional consistency," a weak necessary condition for admissibility discussed in [5].

First we give a sufficient condition for minimaxity.

THEOREM 2.1. *Suppose $(X', Z')'$ has a distribution given by* (1.4). *Then $\hat{\theta}_\phi$ given by* (2.1) *is minimax under $L_j$ if $\phi'(w) \geq 0$ and*

$$0 \leq \phi(w) \leq 2(n+2)^{-1}\left(\frac{\sum\{d_i^{1-j}/c_i\}}{\max\{d_i^{1-j}/c_i\}} - 2\right).$$

PROOF. See the Appendix.

Next we develop a class of generalized Bayes estimators under normality. Suppose the distribution of $(X', Z')'$ is normal with covariance matrix $\sigma^2 \text{diag}(d_1, \ldots, d_p, 1, \ldots, 1)$ and mean vector $(\theta', 0')'$. Consider the following generalized prior distribution:

$$\begin{aligned}
\theta|\lambda, \eta &\sim N_p(0, \eta^{-1}D(\lambda^{-1}C - I)) &&\text{for } \eta = \sigma^{-2}, \\
\lambda &\propto \lambda^a(1-\gamma\lambda)^b I_{[0,1/\gamma]} &&\text{for } \gamma \geq 1, \eta \propto \eta^e.
\end{aligned} \tag{2.2}$$

This is a generalization of priors considered in [3, 4, 10, 14, 18]. The marginal density of $X$, $S$, $\lambda$ and $\eta$ is proportional to

$$\begin{aligned}
&\int \exp\left(-\frac{\eta}{2}\sum\left\{\frac{(x_i - \theta_i)^2}{d_i} + \frac{\lambda}{c_i - \lambda}\frac{\theta_i^2}{d_i}\right\} - \frac{\eta s}{2}\right)\eta^{p+n/2+e}\lambda^{p/2+a} \\
&\quad \times \prod(c_i - \lambda)^{-1/2}(1-\gamma\lambda)^b\, d\theta \\
&\propto \exp\left(-\frac{\eta s}{2}(1+\lambda w)\right)\eta^{p/2+n/2+e}\lambda^{p/2+a}(1-\gamma\lambda)^b,
\end{aligned} \tag{2.3}$$

where $w = x'C^{-1}D^{-1}x/s$. Under the loss $L_j$, the generalized Bayes estimator is given by $E(\eta\theta|X,S)/E(\eta|X,S)$, which can be written, using (2.3),

$$\hat\theta_{\mathrm{GB}} = \left(I - \frac{E(\lambda\eta|X,S)}{E(\eta|X,S)}C^{-1}\right)X = \left(I - \frac{\phi_{\mathrm{GB}}(W)}{W}C^{-1}\right)X.$$

When $p/2 + n/2 + e + 2 > 0$,

(2.4)
$$\int_0^\infty \eta^{p/2+n/2+e+1}\exp\left(-\frac{\eta}{2}\lambda\sum\frac{x_i^2}{c_i d_i} - \frac{\eta s}{2}\right)d\eta$$
$$\propto (1+\lambda w)^{-p/2-n/2-e-2},$$

and we have

(2.5)  $$\phi_{\mathrm{GB}}(w) = \frac{w}{\gamma}\frac{\int_0^1 t^{p/2+a+1}(1-t)^b(1+wt/\gamma)^{-p/2-n/2-e-2}\,dt}{\int_0^1 t^{p/2+a}(1-t)^b(1+wt/\gamma)^{-p/2-n/2-e-2}\,dt},$$

which is well defined for $a > -p/2 - 1$ and $b > -1$. Using an identity which is given by the change of variables $t = (1+w)\lambda/(1+w\lambda)$,

$$\int_0^1 \lambda^\alpha(1-\lambda)^\beta(1+w\lambda)^{-\gamma}\,d\lambda$$
$$= \frac{1}{(w+1)^{\alpha+1}}\int_0^1 t^\alpha(1-t)^\beta\left\{1 - \frac{tw}{w+1}\right\}^{-\alpha-\beta+\gamma-2}\,dt,$$

we have

(2.6)  $$\phi_{\mathrm{GB}}(w) = \frac{w}{\gamma+w}\cdot\frac{\int_0^1 t^{p/2+a+1}(1-t)^b\{1-tw/(w+\gamma)\}^{n/2+e-a-1-b}\,dt}{\int_0^1 t^{p/2+a}(1-t)^b\{1-tw/(w+\gamma)\}^{n/2+e-a-b}\,dt}.$$

The following lemma gives some useful properties of $\phi_{\mathrm{GB}}(w)$.

LEMMA 2.2.  *If $b \geq 0$, $e > -p/2 - n/2 - 2$ and $-p/2 - 1 < a < n/2 + e$, we have for $\phi(w) = \phi_{\mathrm{GB}}(w)$ given by (2.6):*

(i)  *$\phi(w)$ is monotone increasing in $w$.*
(ii) *$\phi(w)/w$ is monotone decreasing in $w$.*
(iii) *$\lim_{w\to\infty}\phi(w) = (p/2+a+1)/(n/2+e-a)$.*
(iv) *$\lim_{w\to 0}\{\phi(w)/w\} = (p/2+a+1)/\{\gamma(p/2+a+b+2)\}$.*

PROOF.  The proof of (i) and (ii) is straightforward using monotone likelihood ratio properties of the densities implied in (2.5) and (2.6). The proof of (iii) and (iv) follows from (2.6) and (2.5), respectively. □

By Lemma 2.2, parts (i) and (iii), and Theorem 2.1, we have immediately the following result.

THEOREM 2.3. *If $b \geq 0$, $e > -p/2 - n/2 - 2$ and $-p/2 - 1 < a < n/2 + e$, then $\hat{\theta}_{\mathrm{GB}}$ is minimax under $L_j$, provided $c_1, \ldots, c_p$ are chosen so that*

$$0 \leq \frac{p/2 + a + 1}{n/2 + e - a} \leq \frac{2}{n+2}\left(\frac{\sum\{d_i^{1-j}/c_i\}}{\max\{d_i^{1-j}/c_i\}} - 2\right).$$

Note if we choose $c_i = d_i/d_p$ under $L_0$, the bound on the RHS is $2(p-2)/(n+2)$. The choices of $a = -2$ and $e = -1$ give a value of $(p-2)/(n+2)$ for the LHS and, hence, for $p \geq 3$ and $n \geq 1$, these choices of $a$ and $e$ give minimax generalized Bayes estimators for any $b \geq 0$ and $\gamma \geq 1$. As Casella [7, 8] indicated, this choice of $c_i$ may be poor from the point of view of the numeric stability of the estimator. It is important to note at this stage that there is substantial flexibility in the choice of $C$ and this flexibility is the key to finding minimax estimators with increased numerical stability. We consider this point further in Sections 3 and 4.

2.1. *A class of simple generalized Bayes minimax estimators.* When $b = n/2 - a + e - 1$ in (2.6), the expression for $\phi_{\mathrm{GB}}(w)$ takes a particularly simple form. In this case,

$$\begin{aligned}
\phi_{\mathrm{GB}}(w) &= \frac{w}{w+\gamma} B(p/2 + a + 2, b + 1) \\
&\quad \times \{B(p/2 + a + 1, b + 1) \\
&\quad\quad - \{w/(w+\gamma)\}B(p/2 + a + 2, b + 1)\}^{-1} \\
&= \frac{\alpha w}{\gamma(\alpha + 1) + w}[= \phi_{\mathrm{SB}}(w), \text{say}],
\end{aligned} \tag{2.7}$$

where $\alpha = (p/2 + a + 1)/(b+1) = (p/2 + a + 1)/(n/2 + e - a)$.

Therefore, our simple generalized Bayes estimator is

$$\hat{\theta}_{\mathrm{SB}} = \left(I - \frac{\alpha}{\gamma(\alpha + 1) + W}C^{-1}\right)X. \tag{2.8}$$

Since $\phi_{\mathrm{SB}}(w)$ is increasing in $w$ and approaches $\alpha$ as $w \to \infty$, we have the following corollary which follows immediately from Theorem 2.1.

COROLLARY 2.4. *$\hat{\theta}_{\mathrm{SB}}$ given by (2.8) is minimax under $L_j$, provided $c_1, \ldots, c_p$ are chosen so that*

$$0 < \alpha \leq \frac{2}{n+2}\left(\frac{\sum\{d_i^{1-j}/c_i\}}{\max\{d_i^{1-j}/c_i\}} - 2\right).$$

In Section 4 we will show that we can always choose $\alpha$, $\gamma$ and $c_1, \ldots, c_p$ to simultaneously achieve minimaxity and an increase in the numerical stability of the least squares estimator.

It is interesting to note that, when $C = D = I_p$, our simple estimator has the form

$$\hat{\theta}_{\mathrm{SB}} = \left(1 - \frac{\alpha}{\gamma(\alpha+1) + X'X/S}\right)X.$$

This is very closely related to Stein's [17] initial class of estimators. He suggested that, for $X \sim N(\theta, I_p)$ with $p \geq 3$, there exist estimators dominating the usual estimator $X$ among a class of estimators of the form $\delta_{a,b} = (1 - b/(a + X'X))X$ for large $a$ and small $b$. Hence, our estimators may be regarded as a variant for the unknown variance case.

Following Stein [17], James and Stein [12] showed that $\delta_{a,b}$ for $a = 0$ and $0 < b < 2(p-2)$ dominates $X$. Since Strawderman [18] derived Bayes minimax estimators, many authors have proposed various minimax (generalized) Bayes estimators. However, the form of these estimators is invariably complicated like our expression (2.6) above. Simple estimators $\delta_{a,b}$ have received little attention although $\delta_{a,b}$, for $a > 0$ and $0 < b < 2(p-2)$, is easily shown to be minimax by using Baranchik's [1] condition. It seems that most statisticians have believed that generalized Bayes estimators which improve on $X$ must have a quite complicated structure. Our result above indicates that this is not so and that generalized Bayes minimax estimators improving on $X$ may indeed have a very simple form.

2.2. *Generalized Bayes estimators for spherically symmetric distributions.* It seems useful to show that the above generalized Bayes results can be extended to the general spherically symmetric case (1.4) in certain situations. What is remarkable about the results is that the resulting generalized Bayes estimators are independent of the form of $f(\cdot)$ and are, hence, identical to those in the normal case. In particular, assume that $C = I$, $\gamma = 1$ and $b = -a - 2$ in the prior given by (2.2). Then the joint density of $\theta$ and $\eta$ is $(\theta' D^{-1} \theta)^{-p/2-a-1} \eta^{-a-1+e}$ because

$$(2.9) \quad \int_0^1 \exp\left(-\frac{\eta\lambda}{2(1-\lambda)}\theta' D^{-1}\theta\right)\eta^{p/2}\left(\frac{\lambda}{1-\lambda}\right)^{p/2}\lambda^a(1-\lambda)^b\,d\lambda$$
$$\propto (\theta' D^{-1}\theta)^{-p/2-a-1}\eta^{-a-1},$$

if $p/2 + a + 1 > 0$. Under quadratic loss $\eta(d-\theta)'(d-\theta)$, the generalized Bayes estimator is given by $E(\eta\theta|X,S)/E(\eta|X,S)$ and we have the generalized Bayes estimator, with respect to our prior,

$$\frac{\int_{R^p}\int_0^\infty \theta\eta^{(n+p)/2-a+e}f(\eta\{X'D^{-1}X+S\})(\theta'D^{-1}\theta)^{-p/2-a-1}\,d\eta\,d\theta}{\int_{R^p}\int_0^\infty \eta^{(n+p)/2-a+e}f(\eta\{X'D^{-1}X+S\})(\theta'D^{-1}\theta)^{-p/2-a-1}\,d\eta\,d\theta}$$
$$= \int_{R^p} \theta(X'D^{-1}X+S)^{-(n+p)/2+a-e-1}(\theta'D^{-1}\theta)^{-p/2-a-1}\,d\theta$$

$$\times \int_0^\infty \eta^{(n+p)/2-a+e} f(\eta) \, d\eta$$

$$\times \left[ \int_{R^p} (X'D^{-1}X + S)^{-(n+p)/2+a-e-1} (\theta' D^{-1}\theta)^{-p/2-a-1} \, d\theta \right.$$

$$\left. \times \int_0^\infty \eta^{(n+p)/2-a+e} f(\eta) \, d\eta \right]^{-1}$$

if $\int_0^\infty \eta^{(n+p)/2-a+e} f(\eta) \, d\eta < \infty$. Note that this does not depend on $f$ and, hence, is equal to the generalized Bayes estimator in the normal case. In the normal case, as seen in Section 2.1, the estimator is well defined if $a > -p/2 - 1$, $b > -1$ and $e > -p/2 - n/2 - 2$. Since $a = -b - 2$, the inequality $-p/2 - 1 < a < -1$ should also be satisfied.

If $b = n/2 - a + e - 1$, which implies $e = -n/2 - 1$, we have a simple generalized Bayes estimator

$$\hat{\theta}_{\text{SB}} = (1 - \alpha/(\alpha + 1 + W))X,$$

where $\alpha = (p/2 + a + 1)/(-a - 1)$ and $W = X'D^{-1}X/S$. Note that $\alpha = (p/2 + a + 1)/(-a - 1)$ can take any positive value because $-p/2 - 1 < a < -1$.

REMARK. The most important point is that, when $a = -b - 2$, $\theta$ and $\eta$ are able to be separated as in (2.9). Furthermore, if $\gamma = 1$, $C = I$ and $a = -b - 2$ are simultaneously not satisfied, the density cannot be so separated. The results in this section are closely related to those in [15].

**3. Condition numbers and numerical stability.** As in Casella [8] and other papers, we use the condition number to measure numerical stability of our ridge-type estimators. This discussion focuses on the stability of estimators of $\beta$ (as opposed to estimators of $\theta$). Recall that our estimators of $\theta$ may be represented as $\hat{\theta}_\phi = (I - tC^{-1})X$, where $t = \phi(w)/w$ and $w = x'C^{-1}D^{-1}x/s$. The vector of regression parameters, $\beta$, is related to the mean vector $\theta$ through the orthogonal matrix $P$ ($\theta = P'\beta$), and the observation vector $X$ in Section 2 is related to the least squares estimator, $\hat{\beta}$, through $X = P'\hat{\beta}$. In this section we study the numerical stability of ridge-type estimators of $\hat{\beta}_\phi$, arising from our improved estimators $\hat{\theta}_\phi$ of $\theta$ through

$$\begin{aligned}
\hat{\beta}_\phi &= P\hat{\theta}_\phi = P(\text{diag}\{d_i^{-1}(1 - t/c_i)^{-1}\})^{-1} P'A'y \\
&= G^{-1}A'y.
\end{aligned} \tag{3.1}$$

By (3.1) $\hat{\beta}_\phi$ may be regarded as a generalized ridge regression estimator $\hat{\beta}_R(K)$ given by (1.5) when we put $k_i = t/\{d_i(c_i - t)\}$.

The condition number of a matrix $H$ is defined by $\kappa(H) = \|H\| \|H^{-1}\|$, where $\|H\| = \sup_{x'x=1} (x'H'Hx)^{1/2} = \max \lambda_i$, where $\lambda_i$ are the eigenvalues

of the positive-definite matrix $H'H$. It follows that if $H$ is a positive-definite matrix, $\kappa(H) = \kappa(H^{-1})$. As indicated in [8] (see also [2]), the condition number measures the numerical sensitivity of the solution of a linear equation $\hat{\beta} = H^{-1}A'y$. In particular, if $\delta\hat{\beta}$ and $\delta(A'y)$ indicate perturbations in $\hat{\beta}$ and $A'y$, respectively,

$$|\delta\hat{\beta}|/|\hat{\beta}| \leq \kappa(H)(|\delta A'y|/|A'y|),$$

where $|\cdot|$ denotes the usual Euclidean norm. For simplicity of notation, we define the condition number of an estimator of the form (3.1) $\kappa(\hat{\beta}_\phi)$ to be equal to the condition number of the matrix $G^{-1}$, $\kappa(G^{-1}) = \kappa(G)$, that is, $\kappa(\hat{\beta}_\phi) = \kappa(G)$.

It follows immediately from the definition of $\kappa(G)$ that (we assume $t < 1$, $c_i \geq 1$)

(3.2) $$\kappa(\hat{\beta}) = d_1/d_p$$

and

(3.3) $$\kappa(\hat{\beta}_\phi) = \frac{\max d_i(1 - t/c_i)}{\min d_i(1 - t/c_i)}.$$

In terms of numerical stability, a smaller condition number implies greater stability. Of course, the condition number given in (3.3) depends on $t = \phi(w)/w$ and, in particular, when $t = 0$, (3.3) reduces to (3.2). We will be interested in finding conditions on the estimator $\hat{\beta}_\phi$ so that, for all possible values of $w$, we have the inequality $\kappa(\hat{\beta}_\phi) \leq \kappa(\hat{\beta})$.

The following result allows condition number improving estimators under two different conditions on $c_1, \ldots, c_p$.

THEOREM 3.1. *Suppose $0 \leq \phi(w)/w \leq t_0 < 1$ for any $w$. Then $\kappa(\hat{\beta}_\phi) \leq \kappa(\hat{\beta})$ for any $w$ if either:*

(i) $c_1 \leq c_2 \leq \cdots \leq c_p$ and

(3.4) $$t_0 \leq \min_{i>j}\left(\frac{c_i c_j(d_1 d_j - d_i d_p)}{c_i d_1 d_j - c_j d_i d_p}\right),$$

*or*

(ii) $c_p > c_1 \geq c_2 \geq \cdots \geq c_{p-1}$ and

(3.5) $$t_0 \leq \min\left(\frac{c_1 c_{p-1}(d_{p-1} - d_p)}{c_1 d_{p-1} - c_{p-1}d_p}, \frac{c_{p-1}c_p(d_1 d_{p-1} - d_p^2)}{c_p d_1 d_{p-1} - c_{p-1}d_p^2}\right).$$

PROOF. See the Appendix.

In the next section we will see that the two conditions above allow us to choose minimax generalized Bayes estimators with increased numerical stability for all full rank designs in the normal case.

**4. Minimaxity and stability.** In this section we show that the results of the previous two sections can be combined to give minimax estimators which simultaneously reduce the condition number relative to the least squares estimator. Then we give a corollary for the simple generalized Bayes estimator $\hat{\theta}_{\text{SB}}$ given by (2.8) under the normal case, because it seems to have practical utility for the general spherically symmetric case. Finally, we add some comments for the case of more general quadratic loss than $L_j$ given by (1.3).

Note that it seems generally desirable to have $c_1 \leq \cdots \leq c_p$ since this implies that the components of $X$ with larger variances get shrunk more. See [8] for an expanded discussion of this point.

Our first result below shows that we may find a minimax condition number improving estimator satisfying $c_1 \leq \cdots \leq c_p$ whenever $\sum \{d_i/d_1\}^{1-j} - 2 > 0$. Note that, when $j \geq 1$, $\sum \{d_i/d_1\}^{1-j} - 2$ is always positive.

THEOREM 4.1. *Suppose $p \geq 3$ and $\sum \{d_i/d_1\}^{1-j} - 2 > 0$. If $d_1 > d_2$, let $\eta_*$ be the unique root such that $\sum \{d_i/d_1\}^\eta = 2$ and let $\eta_{**}$ be any value in $(\max\{0, 1-j\}, \eta_*)$. If $d_1 = d_2$, let $\eta_{**}$ be any value $> \max(1-j, 0)$. Then if $c_i = (d_1/d_i)^{j-1+\eta_{**}}$,*

$$u_+ = 2(n+2)^{-1}\Big(\sum \{d_i/d_1\}^{\eta_{**}} - 2\Big)$$

*and*

$$(4.1) \qquad v_+ = \min_{i>j}\bigg(\frac{c_i c_j (d_1 d_j - d_i d_p)}{c_i d_1 d_j - c_j d_i d_p}\bigg),$$

*the estimator $\hat{\theta}_\phi$ where $0 \leq \phi(w)/w \leq v_+$ for any $w$, $\phi(w)$ is increasing and $\lim_{w \to \infty} \phi(w) \leq u_+$, is minimax under $L_j$ and condition number decreasing, further $c_1 \leq \cdots \leq c_p$.*

PROOF. Since $(d_i/d_1)^\eta$ is strictly decreasing in $\eta$ if $d_i/d_1 < 1$, there exists exactly one root $\eta_*$ of $\sum (d_i/d_1)^\eta = 2$ if $d_2/d_1 < 1$, and that root is strictly larger than $1 - j$. If $d_1 = d_2$, $\sum (d_i/d_1)^\eta > 2$ for any $\eta > 0$. Hence, $\eta_{**} > 1 - j$ and $c_i = (d_1/d_i)^{j-1+\eta_{**}}$ is monotone nondecreasing in $i$. Also from Theorem 2.1 we have minimaxity, provided

$$0 < \phi(w) \leq \frac{2}{n+2}\bigg(\frac{\sum\{d_i^{1-j}/c_i\}}{\max\{d_i^{1-j}/c_i\}} - 2\bigg)$$

$$= \frac{2}{n+2}\Big(\sum\{d_i/d_1\}^{\eta_{**}} - 2\Big) = u_+ \ (> 0).$$

Also by Theorem 3.1(i), since $c_1 \leq c_2 \leq \cdots \leq c_p$, the estimator $\hat{\theta}_\phi$ will have reduced condition number, provided $0 \leq \phi(w)/w \leq v_+$ for any $w$. $\square$

From Theorem 4.1 we easily see the robustness of minimaxity with respect to loss function.

COROLLARY 4.2. *A minimax estimator under $L_j$ for fixed $j$, which is given by Theorem 4.1, retains minimaxity under $L_k$ for $j < k < j + \eta_{**}$.*

For example, suppose $\sum \{d_i/d_1\}^2 - 2 > 0$ and $L_0$ is used. In Theorem 4.1, we can choose $\eta_{**}$ as strictly greater value than 2. Hence, a minimax estimator using such $\eta_{**}$ under $L_0$ retains minimaxity under $L_1$ and $L_2$.

There remains the case where $\sum \{d_i/d_1\}^{1-j} - 2 \leq 0$. Recall that $\sum \{d_i/d_1\}^{1-j} - 2$ is always positive for $j \geq 1$. This case corresponds to the case where there is no spherically symmetric estimator ($c_1 = c_2 = \cdots = c_p$) and, therefore, no estimator with $c_1 \leq c_2 \leq \cdots \leq c_p$ can be minimax (e.g., see [6]). Our solution, while less pleasing in a sense than Theorem 4.1, nevertheless, allows a minimax estimator which reduces the condition number and, hence, increases the stability.

THEOREM 4.3. *Suppose $p \geq 3$ and $\sum \{d_i/d_1\}^{1-j} - 2 \leq 0$. If $p \geq 4$, let $\nu_* \in (0, 1-j)$ be the unique solution of $\sum_{i=1}^{p-1} \{d_i/d_1\}^{\nu} = 2$. Let $\nu_{**}$ be any value in $[0, \nu_*)$. If $p = 3$, choose $\nu_{**} = 0$. Then if $c_i = (d_i/d_{p-1})^{1-j-\nu_{**}}$ for $i = 1, 2, \ldots, p-1$ and $c_p > c_1$,*

$$u_- = 2(n+2)^{-1} \left( \sum_{i=1}^{p-1} \{d_i/d_1\}^{\nu_{**}} - 2 + \frac{c_1}{c_p} \left\{ \frac{d_p}{d_1} \right\}^{1-j} \right)$$

*and*

$$v_- = \min\left( \frac{c_1 c_{p-1}(d_{p-1} - d_p)}{c_1 d_{p-1} - c_{p-1} d_p}, \frac{c_{p-1} c_p (d_1 d_{p-1} - d_p^2)}{c_p d_1 d_{p-1} - c_{p-1} d_p^2} \right),$$

*the estimator $\hat{\theta}_\phi$ where $0 \leq \phi(w)/w \leq v_-$ for any $w$, $\phi(w)$ is increasing and $\lim_{w \to \infty} \phi(w) \leq u_-$, is minimax under $L_j$ and condition number decreasing.*

PROOF. It is easy to see, as in Theorem 4.1, that $\nu_*$, $\nu_{**}$ can be chosen as indicated. In this case, Theorem 2.1 implies minimaxity, provided

$$0 < \phi(w) \leq \frac{2}{n+2} \left( \frac{\sum \{d_i^{1-j}/c_i\}}{\max\{d_i^{1-j}/c_i\}} - 2 \right)$$

$$= \frac{2}{n+2} \left( \sum_{i=1}^{p-1} \{d_i/d_1\}^{\nu_{**}} - 2 + \frac{c_1}{c_p} \left\{ \frac{d_p}{d_1} \right\}^{1-j} \right) = u_- \ (> 0).$$

Theorem 3.1(ii) then implies, since $c_p > c_1 \geq c_2 \geq \cdots \geq c_{p-1}$, that our estimator is condition improving if $0 \leq \phi(w)/w \leq v_-$ for any $w$. □

Combining Lemma 2.2 and Theorems 4.1 and 4.3, we see that versions of Theorems 4.1 and 4.3 are valid for the broad class of generalized Bayes minimax estimators of Theorem 2.3. We omit the straightforward details. We give explicitly a corollary of Theorems 4.1 and 4.3 for our simple generalized Bayes estimator $\hat{\theta}_{\text{SB}}$, because this seems to have practical utility in the general spherically symmetric case.

COROLLARY 4.4. $\hat{\theta}_{\text{SB}} = (I - \alpha/\{\gamma(\alpha+1) + W\}C^{-1})X$ is minimax under $L_j$ and condition number decreasing (and generalized Bayes under normality) if either:

(i) under the setting of Theorem 4.1, $\alpha \leq u_+$ and $\gamma \geq \alpha/\{(\alpha+1)v_+\}$, or
(ii) under the setting of Theorem 4.3, $\alpha \leq u_-$ and $\gamma \geq \alpha/\{(\alpha+1)v_-\}$.

Hence, in the normal case, for any full rank design, we may choose a simple generalized Bayes minimax estimator with increased numerical stability over the least squares estimator $\hat{\beta}$. Further, these estimators remain minimax for all spherically symmetric error distributions.

Finally, we briefly consider the case for general quadratic loss functions $L_Q = \sigma^{-2}(b-\beta)'Q(b-\beta)$ for a positive definite matrix $Q$. Recall that we have assumed $Q = (A'A)^j$ throughout the paper. It is essential for the derivation of minimax estimators with numerical stability in Theorems 4.1 and 4.3 that $A'A$ and $(A'A)^j$ have common eigenvectors. For a general $Q$ which does not have common eigenvectors with $A'A$, let $M$ be a nonsingular matrix which simultaneously diagonalizes $A'A$ and $Q$, where $M$ satisfies

$$M(A'A)^{-1}M' = G = \text{diag}(g_1,\ldots,g_p), \qquad MM' = Q.$$

Let $X = M'\hat{\beta}$ and $\theta = M'\beta$. As in Section 1, we see that $(X', Z')'$ has the joint density $\prod g_i^{-1/2}\sigma^{-N}f(\{(x-\theta)'G^{-1}(x-\theta) + z'z\}/\sigma^2)$ and that the estimation problem of $\theta$ under the squared loss function $(\delta-\theta)'(\delta-\theta)$ is derived as the equivalent canonical problem. Hence, we can have the same minimaxity result in Theorem 2.1 for the shrinkage estimator of the form (2.1) if $g_i$ is substituted for $d_i$. The corresponding generalized ridge estimator becomes

$$(A'A + MKM')^{-1}A'y,$$

where $K = \text{diag}(k_1,\ldots,k_p)$ for $k_i = t/\{g_i(c_i - t)\}$. The eigenvalues of $A'A + MKM'$ (and hence, the condition number), however, cannot be expressed explicitly while, in Section 3, the eigenvalues of $A'A + PKP'$ and the condition number can. As a result, we cannot explicitly construct minimax estimators with numerical stability as in Theorems 4.1 and 4.3.

**5. Numerical results.** In this section we investigate numerically the risk-performance and condition number-performance of our simple Bayes estimator $\hat{\theta}_{\mathrm{SB}}$, given by Corollary 4.4 under $L_j$ for $j = 0, 1, 2$. In the setting of Theorem 4.1, $\eta_{**} = (\max\{0, 1 - j\} + \eta_*)/2$, $\alpha = u_+$ and $\gamma = \alpha/\{(\alpha+1)v_+\}$ are chosen. In the setting of Theorem 4.3, $\nu_{**} = \nu_*/2$, $\alpha = u_-$ and $\gamma = \alpha/\{(\alpha+1)v_-\}$ are chosen. Simulation experiments are done in the following case: $p = 9$, $n = 10$, $D = \mathrm{diag}(\mu^4, \mu^3, \mu^2, \mu, 1, \mu^{-1}, \mu^{-2}, \mu^{-3}, \mu^{-4})$, where $\mu = 1.2, 1.6, 2.0, 2.4$ and $\theta_i = 0, 0.5, 1, 1.5, 2$ for any $i$. The corresponding condition numbers of $D$, $\mu^8$ (and hence, the condition numbers of $\hat{\beta}$ in the original problem), are approximately 4.3, 43, 256 and 1100, respectively. For only two cases, $\mu = 2$ and 2.4 under $L_0$, Theorem 4.3 is applied.

Table 1 shows the relative performance of our simple estimator with respect to risk and expected condition number (ECN), that is:

- $R(\theta, \hat{\theta}_{\mathrm{SB}})/R(\theta, X)$,
- (expected condition number of $\hat{\theta}_{\mathrm{SB}})/\mu^8$,

from 50,000 replications, in each of the above cases. We can draw the following conclusions:

(i) When $\mu$ is large and $\sum\{d_i/d_1\} - 2$ is nonpositive, minimax stable estimators using Theorem 4.3 under $L_0$ have little gain both in the risk improvement and in the ECN improvement. From the numerical results, our contribution of Theorem 4.3 may be just theoretical.

(ii) Under $L_1$, minimax stable estimators have reasonable performances of risk and the ECN, regardless of $\mu$.

(iii) Under $L_2$, when $\mu$ is large, there is little to gain in risk improvement, while there is much to gain in ECN improvement. With better choices of $\eta_{**}$, $\alpha$ and $\gamma$, however, we may have more reasonable performances of risk and ECN.

**Appendix** PROOF OF THEOREM 2.1. Let

$$F(x) = \tfrac{1}{2} \int_x^\infty f(t)\,dt$$

and define

$$E^f[h(X, Z)] = \int\int h(x, z)\sigma^{-N} \prod d_i^{-1/2} f\left(\frac{(x-\theta)'D^{-1}(x-\theta)}{\sigma^2} + \frac{z'z}{\sigma^2}\right) dx\,dz$$

and

$$E^F[h(X, Z)] = \int\int h(x, z)\sigma^{-N} \prod d_i^{-1/2} F\left(\frac{(x-\theta)'D^{-1}(x-\theta)}{\sigma^2} + \frac{z'z}{\sigma^2}\right) dx\,dz,$$

TABLE 1
*Relative performance of our simple estimators under $L_j$ for $j = 0, 1, 2$*

| $\mu$ | $\theta_i$ | $L_0$ | | $L_1$ | | $L_2$ | |
|---|---|---|---|---|---|---|---|
| | | risk | ECN | risk | ECN | risk | ECN |
| 1.2 | 0 | 0.78 | 0.761 | 0.785 | 0.643 | 0.674 | 0.495 |
| | 0.5 | 0.809 | 0.791 | 0.808 | 0.669 | 0.703 | 0.511 |
| | 1 | 0.866 | 0.852 | 0.857 | 0.734 | 0.769 | 0.565 |
| | 1.5 | 0.912 | 0.902 | 0.903 | 0.806 | 0.837 | 0.641 |
| | 2 | 0.941 | 0.935 | 0.934 | 0.861 | 0.889 | 0.718 |
| 1.6 | 0 | 0.894 | 0.95 | 0.778 | 0.597 | 0.955 | 0.417 |
| | 0.5 | 0.917 | 0.963 | 0.801 | 0.637 | 0.956 | 0.425 |
| | 1 | 0.95 | 0.981 | 0.848 | 0.723 | 0.961 | 0.449 |
| | 1.5 | 0.967 | 0.989 | 0.891 | 0.803 | 0.967 | 0.487 |
| | 2 | 0.978 | 0.996 | 0.923 | 0.86 | 0.973 | 0.537 |
| 2 | 0 | 0.994 | 0.999 | 0.778 | 0.594 | 0.995 | 0.415 |
| | 0.5 | 0.994 | 0.999 | 0.807 | 0.652 | 0.995 | 0.419 |
| | 1 | 0.995 | 0.999 | 0.861 | 0.757 | 0.995 | 0.43 |
| | 1.5 | 0.995 | 1 | 0.905 | 0.839 | 0.996 | 0.449 |
| | 2 | 0.995 | 1 | 0.934 | 0.89 | 0.996 | 0.475 |
| 2.4 | 0 | 0.994 | 1 | 0.778 | 0.593 | 0.998 | 0.422 |
| | 0.5 | 0.994 | 1 | 0.819 | 0.679 | 0.999 | 0.424 |
| | 1 | 0.994 | 1 | 0.882 | 0.802 | 1 | 0.431 |
| | 1.5 | 0.994 | 1 | 0.925 | 0.879 | 1 | 0.442 |
| | 2 | 0.994 | 1 | 0.95 | 0.921 | 1 | 0.456 |

where $h(x, z)$ is an integrable function. The identities corresponding to the Stein and chi-square identities for the normal distribution,

(A.1) $$E^f[(X_i - \theta_i)h(X, Z)] = d_i \sigma^2 E^F[(\partial/\partial X_i)h(X, Z)],$$

(A.2) $$E^f[Sg(S)] = \sigma^2 E^F E[ng(S) + 2Sg'(S)],$$

where $S = Z'Z$, are useful in our following proof. We use the version derived in [13], but earlier versions appear in [16] and elsewhere.

The risk of $\hat{\theta}_\phi$ is given by

$$R_j(\theta, \sigma^2, \hat{\theta}_\phi)$$
$$= E^f[(\hat{\theta}_\phi - \theta)' D^{-j}(\hat{\theta}_\phi - \theta)/\sigma^2]$$
(A.3) $$= R(\theta, \sigma^2, X) + E^f\left[\frac{S^2}{\sigma^2} \frac{\sum\{X_i^2/(c_i^2 d_i^j)\}}{(\sum\{X_i^2/(c_i d_i)\})^2} \phi^2\left(\frac{\sum\{X_i^2/(c_i d_i)\}}{S}\right)\right]$$
$$- 2E^f\left[\sum \frac{S}{\sigma^2} \frac{X_i}{c_i d_i^j}(X_i - \theta_i) \sum \frac{X_i^2}{c_i d_i} \phi\left(\frac{\sum\{X_i^2/(c_i d_i)\}}{S}\right)\right].$$

Let $W = X'C^{-1}D^{-1}X/S$. For the second term in (A.3), using (A.2), we have

$$E^f\left[\frac{X'C^{-2}D^{-j}X}{(X'C^{-1}D^{-1}X)^2}\frac{S}{\sigma^2}\left\{S\phi^2\left(\frac{X'C^{-1}D^{-1}X}{S}\right)\right\}\right]$$

$$= E^F\left[\frac{X'C^{-2}D^{-j}X}{X'C^{-1}D^{-1}X}\left\{(n+2)\frac{\phi^2(W)}{W} - 4\phi(W)\phi'(W)\right\}\right].$$

For the third term in (A.3), using (A.1), we have

$$\sum E^f\left[\frac{1}{c_i d_i^j \sigma^2}(X_i - \theta_i)X_i\left(\frac{\sum\{X_i^2/(c_i d_i)\}}{S}\right)^{-1}\phi\left(\frac{\sum\{X_i^2/(c_i d_i)\}}{S}\right)\right]$$

$$= E^F\left[\sum \frac{d_i^{1-j}}{c_i}\frac{\phi(W)}{W} + 2\frac{X'C^{-2}D^{-j}X}{S}\left\{\frac{\phi'(W)}{W} - \frac{\phi(W)}{W^2}\right\}\right].$$

Hence, since $\phi'(w) \geq 0$, we have

$$R_j(\theta, \sigma^2, \hat{\theta}_\phi)$$

$$\leq R_j(\theta, \sigma^2, X) + E^F\left[\frac{\phi(W)}{W}\frac{X'C^{-2}D^{-j}X}{X'C^{-1}D^{-1}X}\right.$$

$$\left.\times\left\{(n+2)\phi(W) - 2\sum \frac{d_i^{1-j}}{c_i}\frac{X'C^{-1}D^{-1}X}{X'C^{-2}D^{-j}X} + 4\right\}\right]$$

$$\leq R_j(\theta, \sigma^2, X) + E^F\left[\frac{\phi(W)}{W}\frac{X'C^{-2}X}{X'C^{-1}D^{-1}X}\right.$$

$$\left.\times\left\{(n+2)\phi(W) - 2\frac{\sum\{d_i^{1-j}/c_i\}}{\max\{d_i^{1-j}/c_i\}} + 4\right\}\right]$$

$$\leq R_j(\theta, \sigma^2, X). \qquad \square$$

PROOF OF THEOREM 3.1. If $c_1 \leq c_2 \leq \cdots \leq c_p$, we have

$$\frac{d_i(1 - t/c_i)}{d_j(1 - t/c_j)} \leq \frac{d_i}{d_j} \qquad \text{for } i < j$$

and

$$\max_t \frac{d_i(1 - t/c_i)}{d_j(1 - t/c_j)} \leq \frac{d_i(1 - t_0/c_i)}{d_j(1 - t_0/c_j)} \qquad \text{for } i > j.$$

Hence, if

$$\max_{i>j}\left(\frac{d_i(1 - t_0/c_i)}{d_j(1 - t_0/c_j)}\right) \leq \frac{d_1}{d_p}$$

or, equivalently,
$$t_0 \le \min_{i>j}\left(\frac{c_i c_j (d_1 d_j - d_i d_p)}{c_i d_1 d_j - c_j d_i d_p}\right),$$
we have
$$\max_t \frac{\max_i d_i(1-t/c_i)}{\min_j d_j(1-t/c_j)} \le \frac{d_1}{d_p},$$
which proves part (i).

Suppose $c_p > c_1 \ge c_2 \ge \cdots \ge c_{p-1}$. Then $d_1(1-t/c_1) \ge \cdots \ge d_{p-1}(1-t/c_{p-1})$ and so
$$\max_t \frac{\max_{i=1,\ldots,p-1} d_i(1-t/c_i)}{\min_{j=1,\ldots,p-1} d_j(1-t/d_j)} \le \frac{d_1(1-t_0/c_1)}{d_{p-1}(1-t_0/c_{p-1})}.$$
Also,
$$\max_t \frac{d_1(1-t/c_1)}{d_p(1-t/c_p)} \le \frac{d_1}{d_p}$$
and
$$\max_t \frac{d_p(1-t/c_p)}{d_{p-1}(1-t/c_{p-1})} \le \frac{d_p(1-t_0/c_p)}{d_{p-1}(1-t_0/c_{p-1})}.$$
Hence, if
$$\max\left(\frac{d_1(1-t_0/c_1)}{d_{p-1}(1-t_0/c_p)}, \frac{d_p(1-t_0/c_p)}{d_{p-1}(1-t_0/c_{p-1})}\right) \le \frac{d_1}{d_p}$$
or, equivalently,
$$t_0 \le \min\left(\frac{c_1 c_{p-1}(d_{p-1} - d_p)}{c_1 d_{p-1} - c_{p-1} d_p}, \frac{c_{p-1} c_p(d_1 d_{p-1} - d_p^2)}{c_p d_1 d_{p-1} - c_{p-1} d_p^2}\right),$$
we have
$$\max_t \frac{\max_i d_i(1-t/c_i)}{\min_j d_j(1-t/c_j)} \le \frac{d_1}{d_p},$$
which proves part (ii). $\square$

CENTER FOR SPATIAL INFORMATION SCIENCE
UNIVERSITY OF TOKYO
7-3-1 HONGO, BUNKYO-KU
TOKYO 113-0033
JAPAN
E-MAIL: maruyama@csis.u-tokyo.ac.jp

DEPARTMENT OF STATISTICS
RUTGERS UNIVERSITY
501 HILL CENTER, BUSCH CAMPUS
PISCATAWAY, NEW JERSEY 08854-8018
USA
E-MAIL: straw@stat.rutgers.edu